\documentclass[11pt]{article}
\usepackage{color}
\usepackage{amsfonts,amssymb,a4, graphicx, bm, makeidx}
\usepackage{amsmath}
\usepackage{soul}
\usepackage{mathrsfs}
\usepackage[all]{xy}
\usepackage[latin1]{inputenc}
\usepackage[dvips]{geometry}
\def\theequation{\thesection.\arabic{equation}}

\newcommand{\qed}{\hfill\rule{3mm}{3mm}}
\newtheorem{corollary}{Corollary}[section]
\newtheorem{theorem}[corollary]{Theorem}

\newtheorem{proposition}[corollary]{Proposition}
\newtheorem{definition}[corollary]{Definition}

\makeatletter \@addtoreset{equation}{section} \makeatother
\definecolor{GreenYellow}{cmyk}{0.15,0,0.69,0}
\definecolor{Yellow}{cmyk}{0,0,1,0}
\definecolor{Goldenrod}{cmyk}{0,0.10,0.84,0}
\definecolor{Dandelion}{cmyk}{0,0.29,0.84,0}
\definecolor{Apricot}{cmyk}{0,0.32,0.52,0}
\definecolor{Peach}{cmyk}{0,0.50,0.70,0}
\definecolor{Melon}{cmyk}{0,0.46,0.50,0}
\definecolor{YellowOrange}{cmyk}{0,0.42,1,0}
\definecolor{Orange}{cmyk}{0,0.61,0.87,0}
\definecolor{BurntOrange}{cmyk}{0,0.51,1,0}
\definecolor{Bittersweet}{cmyk}{0,0.75,1,0.24}
\definecolor{RedOrange}{cmyk}{0,0.77,0.87,0}
\definecolor{Mahogany}{cmyk}{0,0.85,0.87,0.35}
\definecolor{Maroon}{cmyk}{0,0.87,0.68,0.32}
\definecolor{BrickRed}{cmyk}{0,0.89,0.94,0.28}
\definecolor{Red}{cmyk}{0,1,1,0}
\definecolor{OrangeRed}{cmyk}{0,1,0.50,0}
\definecolor{RubineRed}{cmyk}{0,1,0.13,0}
\definecolor{WildStrawberry}{cmyk}{0,0.96,0.39,0}
\definecolor{Salmon}{cmyk}{0,0.53,0.38,0}
\definecolor{CarnationPink}{cmyk}{0,0.63,0,0}
\definecolor{Magenta}{cmyk}{0,1,0,0}
\definecolor{VioletRed}{cmyk}{0,0.81,0,0}
\definecolor{Rhodamine}{cmyk}{0,0.82,0,0}
\definecolor{Mulberry}{cmyk}{0.34,0.90,0,0.02}
\definecolor{RedViolet}{cmyk}{0.07,0.90,0,0.34}
\definecolor{Fuchsia}{cmyk}{0.47,0.91,0,0.08}
\definecolor{Lavender}{cmyk}{0,0.48,0,0}
\definecolor{Thistle}{cmyk}{0.12,0.59,0,0}
\definecolor{Orchid}{cmyk}{0.32,0.64,0,0}
\definecolor{DarkOrchid}{cmyk}{0.40,0.80,0.20,0}
\definecolor{Purple}{cmyk}{0.45,0.86,0,0}
\definecolor{Plum}{cmyk}{0.50,1,0,0}
\definecolor{Violet}{cmyk}{0.79,0.88,0,0}
\definecolor{RoyalPurple}{cmyk}{0.75,0.90,0,0}
\definecolor{BlueViolet}{cmyk}{0.86,0.91,0,0.04}
\definecolor{Periwinkle}{cmyk}{0.57,0.55,0,0}
\definecolor{CadetBlue}{cmyk}{0.62,0.57,0.23,0}
\definecolor{CornflowerBlue}{cmyk}{0.65,0.13,0,0}
\definecolor{MidnightBlue}{cmyk}{0.98,0.13,0,0.43}
\definecolor{NavyBlue}{cmyk}{0.94,0.54,0,0}
\definecolor{RoyalBlue}{cmyk}{1,0.50,0,0}
\definecolor{Blue}{cmyk}{1,1,0,0}
\definecolor{Cerulean}{cmyk}{0.94,0.11,0,0}
\definecolor{Cyan}{cmyk}{1,0,0,0}
\definecolor{ProcessBlue}{cmyk}{0.96,0,0,0}
\definecolor{SkyBlue}{cmyk}{0.62,0,0.12,0}
\definecolor{Turquoise}{cmyk}{0.85,0,0.20,0}
\definecolor{TealBlue}{cmyk}{0.86,0,0.34,0.02}
\definecolor{Aquamarine}{cmyk}{0.82,0,0.30,0}
\definecolor{BlueGreen}{cmyk}{0.85,0,0.33,0}
\definecolor{Emerald}{cmyk}{1,0,0.50,0}
\definecolor{JungleGreen}{cmyk}{0.99,0,0.52,0}
\definecolor{SeaGreen}{cmyk}{0.69,0,0.50,0}
\definecolor{Green}{cmyk}{1,0,1,0}
\definecolor{ForestGreen}{cmyk}{0.91,0,0.88,0.12}
\definecolor{PineGreen}{cmyk}{0.92,0,0.59,0.25}
\definecolor{LimeGreen}{cmyk}{0.50,0,1,0}
\definecolor{YellowGreen}{cmyk}{0.44,0,0.74,0}
\definecolor{SpringGreen}{cmyk}{0.26,0,0.76,0}
\definecolor{OliveGreen}{cmyk}{0.64,0,0.95,0.40}
\definecolor{RawSienna}{cmyk}{0,0.72,1,0.45}
\definecolor{Sepia}{cmyk}{0,0.83,1,0.70}
\definecolor{Brown}{cmyk}{0,0.81,1,0.60}
\definecolor{Tan}{cmyk}{0.14,0.42,0.56,0}
\definecolor{Gray}{cmyk}{0,0,0,0.50}
\definecolor{Black}{cmyk}{0,0,0,1}
\definecolor{White}{cmyk}{0,0,0,0}

\begin{document}
\def\theequation{\thesection.\arabic{equation}}

\def\blu{\color{Blue}}
\def\mag{\color{Maroon}}
\def\red{\color{Red}}
\def\green{\color{ForestGreen}}
\def\prob{{\rm Prob}}

\def\reff#1{(\protect\ref{#1})}

\let\a=\alpha \let\b=\beta \let\ch=\chi \let\d=\delta \let\e=\varepsilon
\let\f=\varphi \let\g=\gamma \let\h=\eta    \let\k=\kappa \let\l=\lambda
\let\m=\mu \let\n=\nu \let\o=\omega    \let\p=\pi \let\ph=\varphi
\let\r=\rho \let\s=\sigma \let\t=\tau \let\th=\vartheta
\let\y=\upsilon \let\x=\xi \let\z=\zeta
\let\D=\Delta \let\F=\Phi \let\G=\Gamma \let\L=\Lambda \let\Th=\Theta
\let\O=\Omega \let\P=\Pi \let\Ps=\Psi \let\Si=\Sigma \let\X=\Xi
\let\Y=\Upsilon

\global\newcount\numsec\global\newcount\numfor
\gdef\profonditastruttura{\dp\strutbox}
\def\senondefinito#1{\expandafter\ifx\csname#1\endcsname\relax}
\def\SIA #1,#2,#3 {\senondefinito{#1#2}
\expandafter\xdef\csname #1#2\endcsname{#3} \else
\write16{???? il simbolo #2 e' gia' stato definito !!!!} \fi}
\def\etichetta(#1){(\veroparagrafo.\veraformula)
\SIA e,#1,(\veroparagrafo.\veraformula)
 \global\advance\numfor by 1
 \write16{ EQ \equ(#1) ha simbolo #1 }}
\def\etichettaa(#1){(A\veroparagrafo.\veraformula)
 \SIA e,#1,(A\veroparagrafo.\veraformula)
 \global\advance\numfor by 1\write16{ EQ \equ(#1) ha simbolo #1 }}
\def\BOZZA{\def\alato(##1){
 {\vtop to \profonditastruttura{\baselineskip
 \profonditastruttura\vss
 \rlap{\kern-\hsize\kern-1.2truecm{$\scriptstyle##1$}}}}}}
\def\alato(#1){}
\def\veroparagrafo{\number\numsec}\def\veraformula{\number\numfor}
\def\Eq(#1){\eqno{\etichetta(#1)\alato(#1)}}
\def\eq(#1){\etichetta(#1)\alato(#1)}
\def\Eqa(#1){\eqno{\etichettaa(#1)\alato(#1)}}
\def\eqa(#1){\etichettaa(#1)\alato(#1)}
\def\equ(#1){\senondefinito{e#1}$\clubsuit$#1\else\csname e#1\endcsname\fi}
\let\EQ=\Eq

\def\pp{{\bm p}}\def\pt{{\tilde{\bm p}}}


\def\\{\noindent}
\let\io=\infty
\def\ee{\end{equation}}
\def\be{\begin{equation}}

\def\VU{{\mathbb{V}}}
\def\EE{{\mathbb{E}}}
\def\N{\mathbb{N}}
\def\U{\mathbb{U}}
\def\GI{{\mathbb{G}}}
\def\TT{{\mathbb{T}}}
\def\C{\mathbb{C}}
\def\CC{{\mathcal C}}
\def\KK{{\mathcal K}}
\def\II{{\mathcal I}}
\def\LL{{\cal L}}
\def\RR{{\cal R}}
\def\SS{{\cal S}}
\def\NN{{\cal N}}
\def\HH{{\cal H}}
\def\GG{{\cal G}}
\def\PP{{\cal P}}
\def\AA{{\cal A}}
\def\BB{{\cal B}}
\def\FF{{\cal F}}
\def\v{\vskip.1cm}
\def\vv{\vskip.2cm}
\def\gt{{\tilde\g}}
\def\E{{\mathcal E} }
\def\EI{{\mathbb E} }
\def\I{{\rm I}}
\def\rfp{R^{*}}
\def\rd{R^{^{_{\rm D}}}}
\def\ffp{\varphi^{*}}
\def\ffpt{\widetilde\varphi^{*}}
\def\fd{\varphi^{^{_{\rm D}}}}
\def\fdt{\widetilde\varphi^{^{_{\rm D}}}}
\def\pfp{\Pi^{*}}
\def\pd{\Pi^{^{_{\rm D}}}}
\def\pbfp{\Pi^{*}}
\def\fbfp{{\bm\varphi}^{*}}
\def\fbd{{\bm\varphi}^{^{_{\rm D}}}}
\def\rfpt{{\widetilde R}^{*}}
\def\A{{{\mathcal O}}}
\def\ef{\mathfrak{f}}
\def\Ti{\mathfrak{T}}
\def\Mi{\mathfrak{M}}
\def\rbfm{\mathrm{\mathbf{m}}}

\def\tende#1{\vtop{\ialign{##\crcr\rightarrowfill\crcr
              \noalign{\kern-1pt\nointerlineskip}
              \hskip3.pt${\scriptstyle #1}$\hskip3.pt\crcr}}}
\def\otto{{\kern-1.truept\leftarrow\kern-5.truept\to\kern-1.truept}}
\def\arm{{}}
\font\bigfnt=cmbx10 scaled\magstep1

\newcommand{\card}[1]{\left|#1\right|}
\newcommand{\und}[1]{\underline{#1}}
\def\1{\rlap{\mbox{\small\rm 1}}\kern.15em 1}
\def\ind#1{\1_{\{#1\}}}
\def\bydef{:=}
\def\defby{=:}
\def\buildd#1#2{\mathrel{\mathop{\kern 0pt#1}\limits_{#2}}}
\def\card#1{\left|#1\right|}
\def\proof{\noindent{\bf Proof. }}
\def\qed{ \square}
\def\trp{\mathbb{T}}
\def\trt{\mathcal{T}}
\def\Z{\mathbb{Z}}
\def\be{\begin{equation}}
\def\ee{\end{equation}}
\def\bea{\begin{eqnarray}}
\def\eea{\end{eqnarray}}
\def\kk{{\bf k}}
\def\Ti{\mathfrak{T}}
\def\Mi{\mathfrak{M}}
\def\begn{\begin{aligned}}
\def\egn{\end{aligned}}
\def\ti{{\rm\bf  t}}\def\mi{{\bm \mu}}
\def\Va{{V^a_{\rm h.c.}}}
\def\Re{{\mathbb{R}}}
\def\T{{\mathcal{T}}}
\def\hL{{\L}}
\def\ev{\mathfrak{e}}
\def\obj{{\rm supp}}\def\fa{\FF}
\def\E{{\cal E}}
\def\EA{{E_\A}}
\def\0{\emptyset}
\def\Ni{\overline{\N}}

\title{A remark on  the Whitney Broken Circuit Theorem}

\author{
\\
Paula M. S. Fialho$^1$,  Emanuel Juliano$^1$, Aldo Procacci$^2$. \\
\\
\small{$^1$Departamento de Ci\^encia da Computa\c{c}\~ao UFMG, }
\small{30161-970 - Belo Horizonte - MG
Brazil}\\
\small{$^2$Departamento de Matem\'atica UFMG,}
\small{ 30161-970 - Belo Horizonte - MG
Brazil}\\
}

\maketitle

\begin{abstract}
\\In the  present note  we show, via the connection between  chromatic polynomial and Potts model, that the Whitney Broken circuit theorem is in fact a special case of a more general identity relating the chromatic polynomial of a graph
$\GI=(\VU,\EE)$ to  sums over forests of  $\GI$  associated to some partition scheme in $\GI$.

\end{abstract}

\section{Introduction}\label{sec1}
Along this paper, $\GI=(\VU, \EE)$ represents a fixed simple graph with vertex set $\VU$ and  edge set $\EE$. For any subset $R\subset \VU$, we denote by $\GI|_R$ the restriction of $\GI$ to $R$, i.e. $\GI|_R=(R, \EE|_R)$ where $\EE|_R=\{\{x,y\}\subset R: \{x,y\}\in \EE\}$.
A graph $\GI'=(\VU',\EE')$ is a subgraph of $\GI=(\VU,\EE)$ if $\VU' \subseteq \VU$ and $\EE' \subseteq \EE$. A subgraph $\GI'$
of $\GI$ is called spanning if $\VU'=\VU$. A graph $\GI=(\VU, \EE)$
is {\it connected}
if for any pair $B, C$ of  subsets of $\VU$ such that
$B\cup C =\VU$ and $B\cap C =\emptyset$, there is at least an edge  $e\in \EE$ such
that $e\cap B\neq\emptyset$ and $e\cap C\neq\emptyset$. A tree is a connected graph with no cycles. A subset $R\subset \VU$ is said to be connected if $\GI|_{R}$ is connected; we denote by $\mathrm{C}_\VU$ the set of all connected subsets of $\GI$ with cardinality greater than one, i.e.
\be\label{civu}
\mathrm{C}_\VU=\{R\subset \VU:\;\mbox{$ \GI|_R$ is connected and $|R|\ge 2$}\}.
\ee

\\Given $R\in \mathrm{C}_\VU$, let us represent by $\mathcal{C}_{\GI|_R}$ the set of all connected spanning subgraphs of ${\GI|_R}$ and by
$\mathcal{T}_{\GI|_R}\subset \mathcal{C}_{\GI|_R}$ the set of  all  spanning trees of $\GI|_R$.

\\Given $n\in \N$, we set $[n]=\{1,2,\dots,n\}$ and given a finite set $X$, we let $|X|$ denote its cardinality. A coloring  of
the vertex set $\VU$ of a graph $\GI=(\VU,\EE)$ using $q \in \N$ colors is a function $\k: \VU\to [q]$ and we denote by $\KK_\VU(q)$ the set of all colorings  of  $\mathbb{V}$ with $q$ colors.  Clearly $\card{\KK_\VU(q)}= q^{\card{\VU}}$. A coloring $\k\in \KK_\VU(q)$ is called {\it proper} if for any edge $\{x,y\}\in \EE$ it holds that
$\k(x)\neq\k(y)$. We represent by $\KK^*_\VU(q)$ the set of all proper colorings of $\VU$ with $q$ colors.
The quantity $P_\mathbb{G}(q):= |\KK^*_\VU(q)|$, i.e. the number of proper colorings with $q$ colors of the graph $\mathbb{G}$ is, as a function of $q$, a polynomial known as the \emph{chromatic polynomial} of $\GI$. The classical representation of the chromatic polynomial $P_\mathbb{G}(q)$ of a graph $\GI=(\VU,\EE)$ is as follows.
\be\label{class}
P_\GI(q)=\sum_{E\subset \EE}(-1)^{|E|} q^{k(E)}
\ee
where $k(E)$ is the number of connected components of the graph $(\VU, E)$ (isolated  vertices are considered trivial connected components, therefore are included).

\\An alternative representation of $P_\mathbb{G}(q)$,  known under the name of Whitney's broken circuit theorem, was derived by Whitney \cite{W}. To state this theorem we need to introduce some notations  regarding  forests of a graph $\GI$.

\\A  connected component of a graph is said non trivial if the cardinality of its vertex set is greater than one, i.e., it contains at least one edge.
A forest in $\GI=(\VU,\EE)$ is a subset $F\subset \EE$ such that all non trivial connected components of the graph $(\VU,F)$ are trees.
We denote by $\mathcal{F}_\GI$ the set of
all forests that are subgraphs of $\GI$ and we agree that the empty set belongs to $\mathcal{F}_\GI$, namely $(\VU, \emptyset)$ is seen as a forest in $\GI$.
Fix now an ordering on the edge set $\EE$ of the graph $\GI$. A circuit $C$ of $\GI$ is a sequence of edges
$\{e_1, \dots, e_n\}$ such that $e_n$ and $e_1$ are adjacent, as well as $e_i$ and $e_{i+1}$ for all $i \in [n-1]$. A circuit can also be seen as a subset of $\EE$. A circuit $C=\{ e_1, \dots, e_n\}$ is called simple if  it does not contain any repeated edges or vertices, except for the vertex $v = e_n \cap e_1$.
Given a simple circuit
$C=\{e_1, \dots, e_n\}$, the set $C\setminus e$ is called a broken circuit of $\GI$ if $e=\max\{e_1, \dots, e_n\}$ (where the maximum is in respect to the ordering of $\EE$ previously  fixed). A forest $F\in  \mathcal{F}_\GI$  is said to be {\it broken circuit free} if it does not contain any broken circuit. Let us represent by $\mathcal{F}^{*}_\GI$ the set of all broken circuit free forests of $\GI$, then the Whitney's broken circuit theorem can now be stated as follows.

\begin{theorem}[Whitney \cite{W}]\label{teowhi}
Given a graph  $\GI=(\VU,\EE) $ , we have
\be\label{whit}
P_\GI(q)= q^{|\VU|}\sum_{F\in \mathcal{F}^{*}_\GI} \left(-{1\over q}\right)^{|F|}.
\ee
\end{theorem}

\\Whitney's broken circuit theorem is relevant because  it provides
a simple combinatorial interpretation of the coefficients of the chromatic polynomial.  Many papers are based on this theorem and several generalization (to hypergraphs, matroids, lattices etc.) have been given (see e.g. \cite{DT} and references therein). In particular,
Theorem \ref{teowhi}  has recently been used in \cite{WQY}  to prove that the number of list-colorings of a graph $\GI=(\VU,\EE)$ is not less than  that of ordinary $q$-colorings of $\GI$ provided that $q>1.135|\EE|$. Even more recently, the broken circuit theorem has been used in
 \cite{JPR} to obtain the best bounds  for the zero-free region of the chromatic polynomial. In \cite{JPR} the  authors show
that for any
graph $\GI$ with maximum degree $\D$, the chromatic polynomial $P_\GI(q)$
is free of zeros   outside the complex disk  $|q|<5.94\Delta$.  Their result improves the previous bounds given in \cite{sok01} and \cite{FP2} which
were obtained instead  via  the  connection between $P_\GI(q)$ and the partition function of the antiferromagnetic Potts model at zero temperature. 

\\In this note we  show that the broken circuit theorem, i.e.  Theorem \ref{teowhi}, is in fact  a special case of a more general representation of  $P_\GI(q)$ in terms of sum over suitable  forest of $\GI$.
To state our  theorem, we need to  introduce  a so-called  {\it partition scheme in $\GI$} (see e.g. Sec. 4 of \cite{JPS} and references therein) according to
the following definition.

\begin{definition}\label{defi}
\\Given $R\in \mathrm{C}_\VU$, a partition scheme in $\GI|_R$ is a map ${ \mathrm{m}_{R}}:\mathcal{T}_{\GI|_R}\to\mathcal{C}_{\GI|_R}$
 such that for all $\t\in \mathcal{T}_{\GI|_R}$, it holds that
$\t\subset {\mathrm{m}_R}(\t)$ and
 $\mathcal{C}_{\GI|_R}=\biguplus_{\tau\in \mathcal{T}_{\GI|_R}}[\tau, \mathrm{m}_R(\tau)]$, where $\biguplus$ means disjoint union
 and $[\tau, \mathrm{m}_R(\tau)]=\{g\in \mathcal{C}_{\GI|_R}: \tau\subseteq g\subseteq {\mathrm{m}_R(\tau)}\}$. Having established, for any $R\subset \mathrm{C}_\VU$, a partition scheme $\mathrm{m}_R$  on ${\cal C}_{\GI|_R}$, the family ${\mathrm{\mathbf{m}}}=\{\mathrm{m}_R\}_{R\in \mathrm{C}_\VU}$ is said to be a partition scheme in $\GI$ and, for any non trivial $\t$ subtree of $\GI$, we will denote shortly
$\mathrm{ \mathbf{m}}(\t)=\mathrm{m}_{V_\t}(\t)$ (here $V_\t$ denotes the vertex set of $\t$).
\end{definition}

\\We are now in the position to state the main result of this note.
\begin{theorem}\label{partschem} Let $\GI=(\VU,\EE)$ be a graph and let  $P_\GI(q)$ be the chromatic polynomial of $\GI$. Let $\mathrm{\mathbf{m}}$ be a partition scheme in $\GI$  according to Definition \ref{defi}.
 Then
\be\label{pgen}
P_\GI(q)=q^{|\VU|} \sum_{F\in \mathcal{F}^{\bf m}_\GI} \left(-{1\over q}\right)^{|F|}
\ee
where $\mathcal{F}^{\bf m}_\GI$ is the subset of  $\mathcal{F}_\GI$ formed by the empty set and  by all non trivial forests $F= \{\t_1, \dots ,\t_k\}$ in which  $\mathrm{\mathbf{m}}(\t_i)=\t_i$  for each tree $\t_i$, $i \in [k]$.
\end{theorem}

\\{\bf Remark}. It is noteworthy to observe that identity \reff{pgen} holds independently of the partition scheme ${\bf m}$ in $\GI$ used.
Even more impressive is that despite the fact that the set $\mathcal{F}_\GI^{\mathrm{\mathbf{ m}}}$ depends strongly on the partition scheme $\rbfm$
, i.e. in general $\mathcal{F}_\GI^{\rbfm}\neq \mathcal{F}_\GI^{\rbfm'}$ if $\rbfm \neq \rbfm'$, the coefficients of \ref{pgen}, namely
$|\{F\in \mathcal{F}_\GI^{\rbfm}: |F|=k\}|$, does not depend on $\rbfm$.
Various partition schemes are available (see e.g. \cite{JPS}, Sec. 4 and references therein). In particular, in  statistical mechanics, as far as we know, two partition schemes have been used: the first one is the so called Penrose partition scheme,  originally presented in \cite{pen67} and used in \cite{FP} to improve the Dobrushin criterion \cite{dob96} for the convergence of the cluster expansion on the abstract polymer gas; the second one  is  the {\it minimal-tree partition scheme}, described  e.g. in the proof of Lemma 2.2. in \cite{scosok05} and used in \cite{PY} to improve the lower bound of the convergence radius of the Mayer Series on a gas of
classical particles in the continuum interacting via a stable and regular pair potential. The Penrose partition scheme has also been
used in combinatorics, e.g. in  \cite{JPS} it was a key tool to find zero-free regions in the complex plane at large $|q|$ for the
multivariate Tutte polynomial (see there Proposition 4.4.and Lemma 4.5). Last but not least,
the improvement of  the Lov\'asz local lemma given in  \cite{BFPS} also relies crucially on the Penrose partition scheme.
\vv

\\To conclude the introduction, we show as promised that the broken circuit theorem, i.e. Theorem \ref{teowhi} is in fact
a special case of Theorem \reff{partschem}. Namely, the broken circuit theorem is in fact Theorem \ref{teowhi} when the partition scheme $\rbfm$
is the so-called {\it  minimal-tree partition scheme}.
To see this, let us construct explicitly the above mentioned partition scheme in $\GI$.

\vv

\\{\it Minimal-tree partition scheme.} We start by choosing a total order $\succ$ in the set of edges $\EE$ of $\GI$. Now, for any $R\in \mathrm{C}_\VU$,  let ${\m_R}: \mathcal{T}_{\GI|_R}\to\mathcal{C}_{\GI|_R}$ be the map that associates to the
tree $\tau=(R, E_\t)\in \mathcal{T}_{\GI|_R}$  the graph ${\m_R}( \tau )\in \mathcal{C}_{\GI|_R}$
whose edge set is $E_{{\m_R}(\tau)}= E_\tau \cup E_\tau^{*}$, where $E_\tau^{*}$ is composed by all edges $\{x,y\}\subset \EE|_R\setminus E_\t$
such that $\{x,y\} \succ \{z,u\} $ for every edge $\{z,u\} \in E_\tau$ belonging to the path from $x$ to $y$ through $\tau$.
\vv
\\It is easy to show that ${\m_R}$  is, for any $R\in C_\VU$,  a partition scheme in $\mathcal{C}_{\GI|_R}$  (see e.g.  \cite{PY}, \cite{PYob}, \cite{scosok05}) and thus the family ${\bm \mu}=\{\mu_R\}_{R\in \mathrm{C}_\VU}$ is a partition scheme in $\GI$. Then, according to Theorem \ref{partschem}, we can write
\be\label{pgbcf}
P_\GI(q)= q^{|\VU|}\sum_{F\in \mathcal{F}^{\bm\m}_\GI} \left(-{1\over q}\right)^{|F|},
\ee
where $\bm \mu$ is the minimal-tree partition scheme  and $\mathcal{F}^{\bm\m}_\GI$ is the set of all forests $F$ in $\GI$  whose non trivial trees $\t \in F$ are such that
 ${\bm\m}(\t)=\t$.
Observe now that
 for any non trivial tree $\t\subset\GI$  the condition $\bm\m(\t)=\t$ simply means that for any $\{x,y\}\in \EE|_R\setminus E_\t$,
the path $p^\t_{xy}$ in $\t$ from $x$ to $y$ is such that $\{x,y\}$ is never the largest edge in the set
$\{\{x,y\}\cup \{e\}_{e\in p^\t_{xy}}\}$.
Therefore, the condition ${\bm\m}(\t)=\t$ for any tree in a forest $F\in\mathcal{F}^{\bm\m}_{\GI}$ can be rephrased by saying that $F$ does not contain broken circuits. That is to say,
$ \mathcal{F}^{\bm\m}_\GI= \mathcal{F}^*_\GI$ and the sum  r.h.s. of \reff{pgbcf} coincides term by term with the sum of the r.h.s. of \reff{whit}.

\section{Proof of Theorem \ref{partschem}}\label{secxi}
\\\\It is long known (see \cite{KE} and references therein) that  $P_\mathbb{G}(q)$ also
coincides with the partition function of the antiferromagnetic
Potts model with $q$ states on $\mathbb{G}$ at zero temperature.
The Potts model with $q$ states on a finite graph $\mathbb{G}=(\mathbb{V}, \mathbb{E})$
(a heavily  studied spin system in statistical mechanics) is defined by a set
of random variables (a.k.a. ``spins") $\{\sigma_x\}_{x\in \VU}$. For any vertex $x\in \VU$, the spin  $\s_x$
takes values  in the set of integers
$[q]$ and a {\it configuration} $\bm \s$ of the system is a function $\bm \s: \VU\to [q]$. Interpreting   $[q]$  as a set of $q$ colors, a spin configuration
$\bm \s $ can be seen as an element of $\KK_\VU(q)$.
The {\it energy} of a spin configuration $\bm\s\in \KK_\VU(q)$ is defined as
\begin{equation*}
H_{\GI}(\bm \s) =
-J\sum_{\{x,y\}\in \mathbb{E}}\d_{\sigma_x \sigma_y}
\end{equation*}
where  $J\in \mathbb{R}$ and
$\d_{\sigma_x\sigma_y}$ is the Kronecker symbol which is equal to one when $\sigma_x=\sigma_y$ and
zero otherwise.
When    $J>0$ the model is called
{\it ferromagnetic}, while if  $J< 0$ the model is said to be  {\it antiferromagnetic}.
{\it The probability} to find the system in the configuration $\bm\s$ is defined as
\be\label{1.2.0}
{\rm Prob}(\bm \s)= \frac{e^{-\b H_{\GI}(\bm \s)}}{ Z_{\GI}(q, \b)}
\ee
where $\b\ge 0$ is the inverse temperature. The normalization constant $Z_{\GI}(q, \b)$ in the denominator
is called the  {\it partition function} and is explicitly given by
\begin{equation}\label{zpotts}
Z_{\GI}(q, \b)=\sum_{\bm\s\in \KK_\VU(q)} e ^{- \b H_{\GI}(\bm \s)}.
\end{equation}
The case $J<0$ and $\b=+\infty$ is the {\it antiferromagnetic} and
{\it zero temperature} Potts model with $q$-states. In this case,
when  a product $0\cdot\infty$  appears in the factor $\b H_{\GI}(\bm \s)$,
it must be interpreted as zero.  With this convention, the only configurations $\bm\s$ with non zero probability are  those for which
$\d_{\s_x\s_y}=0$ for all $\{x,y\}\in \EI$, i.e. those configurations in
which adjacent spins have different values (or colors).
In other words, the
antiferromagnetic  Potts model at zero temperature is such that ${\rm Prob}(\bm\s)\neq 0$ if and only if the configuration $\bm\s$ is a proper coloring (i.e. $\bm\s\in \KK^*_\VU(q)$).
Clearly, for any $\bm\s\in \KK^*_\VU(q)$  we have that
$\exp\{-\b H_{\GI}(\bm\s\}|_{\b=+\infty, J<0}=1$, so that, according to (\ref{1.2.0}), all colorings in $\KK^*_\VU(q)$ have the same probability
$1/Z_\GI(q, \b)$,  and thus the normalization constant  $Z_\GI(q, \b)$  coincides with $|\KK^*_\VU(q)|$, i.e. we get that
\be\label{Potts0}
 Z_\GI(q,\b)|_{\b=+\infty, J<0} =P_\GI(q).
\ee
\vskip.2cm
\\The connection between chromatic polynomial and antiferromagnetic Potts model at zero temperature shown by Formula \reff{Potts0}
has been explored in recent times to study the location of the zeros of the $P_\GI(q)$ in the complex plane
(see \cite{sok01,FP2}).

\\Let us enunciate two rather standard results in statistical mechanics. The first  is  basically related to the well known
high-temperature expansion for  the partition function of a general discrete spin system while the second is  a long known identity discovered by  Oliver  Penrose \cite{pen67} during  with the study of
the mathematical properties of continuous particle systems in the high-temperature and  low-density regime.

\begin{proposition}\label{pro1}
let  $Z_\GI(q,\b)|_{\b=+\infty, J<0}$  be the partition function of the antiferromagnetic Potts model with $q$ state at zero temperature. Then
\be\label{zqx}
Z_\GI(q,\b)|_{\b=+\infty, J<0}=q^{|\VU|}\Xi_\GI(q),
\ee
where
\be\label{XiG}
\Xi_\GI(q) =1+\sum_{k\ge 1}\sum_{\{R_1,\dots, R_k\}\atop
 R_i\in \mathrm{C}_\VU, ~R_i\cap R_j=\emptyset} \prod_{i =1}^{k}z(R_i,q)
\ee
is the grand canonical partition function of a ``polymer gas'', in which the polymers  are connected subsets $R\subset \mathbb{V}$, with cardinality $|R|\ge 2$  (i.e. they are the elements of the set
$\mathrm{C}_\VU$ defined in \reff{civu}),  subjected to a non-intersection constraint (hard-core interaction) and with weights (a.k.a. activities) $z(R,q)$ given by
\be\label{actrq}
z(R,q)={1\over q^{|R|-1}} \sum_{g\in\mathcal{C}_{\GI|_R}}(-1)^{|E_g|}
\ee
where $\mathcal{C}_{\GI|_R}$ denotes  the set of all connected spanning subgraphs of $\GI|_R$.
\end{proposition}
A sketch of the proof of Proposition \ref{pro1}  can  be found e.g. in \cite{sok01} and in \cite{JPS}.
However, since this representation of $ Z_\GI(q,\b)|_{\b=+\infty, J<0}$ is crucial for the proof of
Theorem \ref{partschem}, we will provide a detailed proof  of it  in the appendix.

\begin{proposition}
Let ${\rbfm}$  be a   partition scheme in $\GI=(\VU,\EE)$ according to Definition \ref{defi}, and let  $\{w_e\}_{e\in \EE}$ be
any collection of real numbers, then, for any  $R\in \mathrm{C}_\VU$, we have
\be\label{pen}
\sum_{g\in\mathcal{C}_{\GI|_R}} \prod_{e\in E_g} w_e =
\sum_{\t \in \mathcal{T}_{\GI|_R}} \prod_{e\in E_\t} w_e \prod_{e \in E_{{\rbfm}(\t)} \setminus E_\t} (1+w_e),
\ee
where $E_g$, $E_\t$, and $E_{{\rbfm}(\t)}$ represent the edge sets of $g$, $\t$, and ${\rbfm}(\t)$ respectively.
\end{proposition}
The proof of identity \reff{pen} is straightforward and can be seen e.g. in Sec. 4.1 of \cite{FP}.

\\We are now in the position to  prove Theorem \ref{partschem}.
Let us turn initially
our attention to the Penrose identity \reff{pen} and observe that when $w_e=-1$ for all $e\in \EE|_R$,
the l.h.s. of \reff{pen} becomes $\sum_{g\in\mathcal{C}_{\GI|_R}}(-1)^{|E_g|}$, which, recalling the definition of $z(q,R)$ given in \reff{actrq} is exactly the sum in the r.h.s. of \reff{actrq}. On the other hand  in the r.h.s. of the Penrose identity  \reff{pen}, when  $w_e=-1$ for all $e\in \EE|_R$, the last product   kills the contribution of any tree $\t$ such that $ E_{{ \rbfm}(\t)}\setminus E_\t\neq\0$, so that  the only contributions on the sum over trees in $\mathcal{T}_{\GI|_R}$ come from those trees $\t$ such that ${\rbfm}(\t)=\t$. Since
for any $\t\in \mathcal{T}_{\GI|_R}$ we have that $|\t|=|R|-1$, as far as the case $w_e=-1$ for all $e\in \EE|_R$ is concerned, the Penrose identity \reff{pen} can be written as
\be\label{ide}
\sum_{g\in\mathcal{C}_{\GI|_R}}(-1)^{|E_g|}= 
(-1)^{|R|-1} \sum_{\t\in \mathcal{T}_{\GI|_R}\atop {\rbfm}(\t)=\t}1.
\ee
Hence we can rewrite
the activity $z(R,q)$ of a polymer $R$, given in \reff{actrq}, as
\be\label{actsche}
z(R,q)=\sum_{\t \in \mathcal{T}_{\GI|_R}\atop {\rbfm}(\t)=\t}\left(-{1\over q}\right)^{|\t|},
\ee
where ${\rbfm}$ is any partition scheme on $\GI$.
Then, recalling \reff{XiG} we get
\be\label{whigen}
\Xi_\GI(q)= 1+\sum_{k\ge 1}\sum_{\{R_1,\dots, R_k\}\atop
 R_i\in \mathrm{C}_\VU, ~R_i\cap R_j=\emptyset} \prod_{i =1}^{k}\sum_{\t_i \in \mathcal{T}_{\GI|_{R_i}}\atop {\rbfm }(\t_i)=\t_i}\left(-{1\over q}\right)^{|\t|}
=\sum_{F\in \mathcal{F}^{\rbfm}_\GI} \left(-{1\over q}\right)^{|F|},
\ee
\vskip.2cm
\\where $\mathcal{F}^{\rbfm}_\GI$ is the subset of  $\mathcal{F}_\GI$ formed by  the empty set
(contributing  with  the factor ``1" in the sum in the r.h.s. of \reff{whigen})
and  by all non trivial forests $F= \{\t_1, \dots ,\t_k\}$ where each tree $\t_i$, $i \in [k]$, is such that $\rbfm(\t_i)=\t_i$.
Inserting finally \reff{whigen} into  \reff{zqx} and recalling \reff{Potts0}, we get \reff{pgen}.
\vv
\section*{Acknowledgements}
It is a pleasure to  thank Benedetto Scoppola for discussions and  comments. P.M.S.F. and  E.J.  were supported by
FAPEMIG (Funda\c{c}\~ao de Amparo \`a Pesquisa do Estado de Minas Gerais).

\renewcommand{\theequation}{A.\arabic{equation}}
\setcounter{equation}{0}

\section*{Appendix. Proof of Proposition \ref{pro1}}
We start by  rewriting $Z_\GI(q)|_{\b=+\infty, J<0}$ as follows.
\be\label{2.4b}
Z_\GI(q)|_{\b=+\infty, J<0}=\sum_{\bm \s\in \KK_\VU(q)}
\exp \Big\{- \sum_{\{x,y\}\subset \VU}J_{xy}\d_{\sigma_x \sigma_y} \Big\}
\ee
where
\be\label{Jxy}
 J_{xy}=\begin{cases} +\infty &{\rm if}~ \{x,y\}\in \mathbb{E}\\
0 & {\rm otherwise}
\end{cases}
\ee
reminding that a product $0\cdot\infty$ in  the r.h.s. of (\ref{2.4b}) must be interpreted as zero.

\\Using the so-called Mayer trick, we can write
$$
\begin{aligned}
\exp \Big\{- \sum_{\{x,y\}\subset \VU}J_{xy}\d_{\sigma_x \sigma_y} \Big\} & =
\prod_{\{x,y\}\subset  \VU}\Big[(e^{-J_{xy}\d_{\sigma_x \sigma_y}}-1)+1]\\
& = \sum_{g\in \GG_\VU} \prod_{\{x,y\}\in E_g}(e^{-J_{xy}\d_{\sigma_x \sigma_y}}-1)
\\
&=1+ \sum_{k\ge 1}\sum_{\{R_1,\dots, R_k\}\subset \VU\atop
R_i\cap R_j=\emptyset,~|R_i|\geq 2} \prod_{i =1}^{k}\sum_{g\in G_{R_i}}\prod_{\{x,y\}\in E_g}(e^{-J_{xy}\d_{\sigma_x \sigma_y}}-1)\\
&=1+\sum_{k\ge 1}\sum_{\{R_1,\dots, R_k\}\subset \VU\atop
R_i\cap R_j=\emptyset,~|R_i|\geq 2} \prod_{i =1}^{k}\r(R_i,\bm\s_{R_i})
\end{aligned}
$$
where in the second line  we have denoted by $\GG_\VU$ the set of all graphs whose vertex set is $\VU$ and, for any $g\in \GG_\VU$, $E_g$ denotes the  edge set of $g$. In the third line, for any $R\subset \VU$ such that $|R|\ge 2$, $G_R$ denotes  the set of all connected graphs with vertex set $R$. Finally, in the fourth line
we have set, for any $R\subset \VU$ such that $|R|\ge 2$,
\be\label{rors}
\r(R,\bm\s_{R})= \sum_{g\in G_{R}}\prod_{\{x,y\}\in E_g}(e^{-J_{xy}\d_{\sigma_x \sigma_y}}-1)
\ee
where for any $R\subset \VU$, $\bm\s_{R}$ is the restriction of $\bm \s$ to $R$.

\\Recalling \reff{Potts0}, we have that the chromatic polynomial $P_\GI(q)$ can be written as
$$
\begin{aligned}
P_\GI(q)& =\sum_{\bm \s \in \KK_\VU(q)}\Bigg[1+\sum_{k\ge 1}\sum_{\{R_1,\dots, R_k\}\subset \VU\atop
R_i\cap R_j=\emptyset,~|R_i|\geq 2} \prod_{i =1}^{k}\r(R_i,\s_{R_i})\Bigg]\\
& = q^{|\VU|}+\sum_{k\ge 1}\sum_{\{R_1,\dots, R_k\}\subset \VU\atop
R_i\cap R_j=\emptyset,~|R_i|\geq 2} \prod_{i =1}^{k}\sum_{\bm\s_{R_i} \in \KK_{R_i}(q)}\r(R_i,\bm\s_{R_i}) q^{|\VU|-\sum_{i=1}^k|R_i|}
\end{aligned}
$$
According to the last equality above we can rewrite $P_{\GI}(q)$ as follows.
\be\label {Zgiq1}
P_{\GI}(q) = q^{|\VU|} \Bigg(1+\sum_{k\ge 1}\sum_{\{R_1,\dots, R_k\}\subset \VU\atop
R_i\cap R_j=\emptyset,~|R_i|\geq 2} \prod_{i =1}^{k}{1\over q^{|R_i|}}\sum_{\bm\s_{R_i}\in \KK_{R_i}(q)}\r(R_i,\bm\s_{R_i})\Bigg)
\ee

\\Now observe that,  given  $R\subset \VU$ such that $|R|\ge 2$ we have that
\be\label{actind}
\sum_{\bm\s_{R}\in \KK_{R}(q)}\r(R,\bm\s_{R})=
\begin{cases}\displaystyle{q\sum\limits_{g\in \mathcal{C}_{\GI|_R}}}(-1)^{|E_g|} & \mbox{if $\GI|_R$ is connected}
\\
0 & {\rm otherwise}
\end{cases}
\ee
where $\mathcal{C}_{\GI|_R}$ denotes  the set of all connected spanning subgraphs of $\GI|_R$.
Indeed, recalling \reff{rors}, we have
$$
\begin{aligned}\sum_{\bm\s_{R}\in \KK_{R}(q)}\r(R,\bm\s_{R})& =\sum_{\bm\s_{R}\in \KK_{R}(q)}\sum_{g\in G_{R}}\prod_{\{x,y\}\in E_g}(e^{-J_{xy}\d_{\sigma_x \sigma_y}}-1)\\
& = \sum_{\bm\s_{R}\in \KK_{R}(q)}\sum_{g\in G_{R}}\prod_{\{x,y\}\in E_g}\d_{\sigma_x \sigma_y}(e^{-J_{xy}}-1)\\
& =\sum_{g\in G_{R}}\left[\sum\limits_{\bm\s_{R}\in \KK_{R}(q)}\prod\limits_{\{x,y\}\in E_g} \d_{\sigma_x \sigma_y}\right]\prod_{\{x,y\}\in E_g}(e^{-J_{xy}}-1)\\
&= q\sum_{g\in G_{R}}\prod_{\{x,y\}\in E_g}(e^{-J_{xy}}-1)
\end{aligned}
$$
where the  last equality above  holds because, for any $g\in G_R$,
$$
\sum\limits_{\bm\s_{R}\in \KK_{R}(q)}\prod\limits_{\{x,y\}\in E_g} \d_{\sigma_x \sigma_y}=q
$$
Finally, recalling the definition of $J_{xy}$ given in \reff{Jxy}, we have that
$$
\sum_{g\in G_{R}}\prod_{\{x,y\}\in E_g}(e^{-J_{xy}}-1)= \begin{cases}\displaystyle{\sum\limits_{g\in \mathcal{C}_{\GI|_R}}}(-1)^{|E_g|} & \mbox{if $\GI|_R$ is connected}
\\
0 & {\rm otherwise}
\end{cases}
$$
Therefore inserting (\ref{actind}) in (\ref{Zgiq1}) and recalling that
$\mathrm{C}_\VU= \{R\subset \VU: |R|\ge 2~{\rm and}~ \GI|_R ~{\rm is}~{\rm connected}\}$ (i.e. $\mathrm{C}_\VU$ is  the set of all connected subsets of $\VU$  with cardinality greater than one), we get
$$
P_\GI(q)=q^{|\VU|}\Xi_\GI(q)
$$
where
$$
\Xi_\GI(q) =1+\sum_{k\ge 1}\sum_{\{R_1,\dots, R_k\}\atop
 R_i\in \mathrm{C}_\VU, ~R_i\cap R_j=\emptyset} \prod_{i =1}^{k}z(R_i,q)
$$
with
$$
z(R,q)={1\over q^{|R|-1}} \sum_{g\in\mathcal{C}_{\GI|_R}}(-1)^{|E_g|}
$$
This concludes the proof of Proposition \ref{pro1}.

\end{document}